\newcommand{\trainingData}{\mathcal{D}}
\newcommand{\learningRate}{\eta}
\newcommand{\loss}{\mathcal{L}}
\newcommand{\objective}{\mathcal{O}}
\newcommand{\stress}{\sigma}

\newcommand{\crystalPullingRate}{v_\text{c}}
\newcommand{\thermalGradient}{G_\text{z}}
\newcommand{\inducedHeatEM}{\dot{q}^\text{EM}}
\newcommand{\triplePhaseLine}{L}

\newcommand{\radialThermalGradient}{G_\text{r}}
\newcommand{\deflection}{\delta}
\newcommand{\exceedStress}{\stress_\text{ex}}
\newcommand{\voltage}{U}
\newcommand{\EMhomogeneity}{H^\text{EM}}
\newcommand{\VoronkovCriterion}{\Gamma}

\newcommand{\modelOutput}{\mathcal{Y}}
\newcommand{\outputRadialThermalGradient}{\modelOutput_1}
\newcommand{\outputDeflection}{\modelOutput_2}
\newcommand{\outputExceedStress}{\modelOutput_3}
\newcommand{\outputVoltage}{\modelOutput_4}
\newcommand{\outputEMhomogeneity}{\modelOutput_5}
\newcommand{\outputVoronkovCriterion}{\modelOutput_6}

\newcommand{\approptoinn}[2]{\mathrel{\vcenter{
  \offinterlineskip\halign{\hfil$##$\cr
    #1\propto\cr\noalign{\kern2pt}#1\sim\cr\noalign{\kern-2pt}}}}}

\newcommand{\modelInput}{\mathcal{X}}
\newcommand{\inputCrystalRadius}{\modelInput_{1}}
\newcommand{\inputPullingRate}{\modelInput_{2}}
\newcommand{\inputSideSlits}{\modelInput_{3}}
\newcommand{\inputSideSlitLength}{\modelInput_{4}}
\newcommand{\inputSideSlitWidth}{\modelInput_{5}}
\newcommand{\inputMainSlitWidth}{\modelInput_{6}}
\newcommand{\inputBottomAngle}{\modelInput_{7}}
\newcommand{\inputFrequency}{\modelInput_{8}}
\newcommand{\inputReflectorHeight}{\modelInput_{9}}
\newcommand{\inputReflectorRadius}{\modelInput_{10}}
\newcommand{\inputReflectorPosition}{\modelInput_{11}}
\newcommand{\inputReflectorEmissivity}{\modelInput_{12}}

\newcommand{\samplesLHS}{S}
\newcommand{\dimensionsLHS}{n}
\newcommand{\coeffDetermination}{\text{R}^2}

\documentclass{WileyMSP-template}

\usepackage{amssymb}
\usepackage{amsmath}
\usepackage[acronyms]{glossaries}
\usepackage{siunitx}
\usepackage{booktabs}
\usepackage{mathtools}
\usepackage{makecell}
\usepackage{xfrac}
\usepackage{mathrsfs}
\usepackage{interval}
\usepackage{multirow}
\usepackage{algorithm}
\usepackage{algpseudocode}
\usepackage{comment}

\newacronym{FZ}{FZ}{Floating Zone}
\newacronym[plural=FEMs]{FEM}{FEM}{Finite Element Model}
\newacronym{ETP}{ETP}{External Triple Point}
\newacronym{TPL}{TPL}{Three-Phase Line}
\newacronym{LPS}{LPS}{Lateral Photovoltage Scanning}
\newacronym{EM}{EM}{Electromagnetic}
\newacronym{ML}{ML}{Machine Learning}
\newacronym[longplural={Quantities of Interest}, plural=QoIs]{QoI}{QoI}{Quantities of Interest}
\newacronym[plural=NNs]{NN}{NN}{Neural Network}
\newacronym{FFNN}{FFNN}{Feed-Forward Neural Network}
\newacronym[plural=DNNs]{DNN}{DNN}{Deep Neural Network}
\newacronym[plural=DEs]{DE}{DE}{Deep Ensemble}
\newacronym[plural=GAs]{GA}{GA}{Genetic Algorithm}
\newacronym{MSE}{MSE}{Mean Squared Error}
\newacronym{ReLU}{ReLU}{Rectified Linear Activation}
\newacronym{DoE}{DoE}{Design of Experiments}
\newacronym{LHS}{LHS}{Latin Hypercube Sampling}
\newacronym{NSGA}{NSGA}{Non-dominated Sorting Genetic Algorithm}
\newacronym{MOO}{MOO}{Multi-Objective Optimization}
\newacronym{TPE}{TPE}{Tree-structured Parzen Estimator}
\newacronym{SGD}{SGD}{Stochastic Gradient Descent}
\newacronym{IKZ}{IKZ}{\textit{Leibniz-Institut für Kristall\-züchtung}}

\DeclareUnicodeCharacter{2212}{-}

\begin{document}

\pagestyle{fancy}

\title{Data-Driven Multi-Objective Optimization of Large-Diameter Si Floating-Zone Crystal Growth}
\maketitle

\author{Lucas Vieira*}
\author{Milena Petković}
\author{Robert Menzel}
\author{Natasha Dropka*}

\begin{affiliations}
L. Vieira, M. Petković, R. Menzel, N. Dropka\\
Leibniz-Institut für Kristallzüchtung, Max-Born-Straße 2, Berlin, 12489, Germany\\
Email Addresses: lucas.vieira@ikz-berlin.de, milena.petkovic@ikz-berlin.de, robert.menzel@ikz-berlin.de, natascha.dropka@ikz-berlin.de\\
\end{affiliations}


\begin{abstract}
Floating Zone (FZ) silicon crystal growth is essential for high-power electronics and advanced detection systems.
The increasing pressure to scale up the process is challenging due to competing objectives. This study presents a surrogate-based optimization framework to address Multi-Objective Optimization (MOO) in FZ growth, considering eight relevant objectives related to productivity, geometrical and growth parameters, and crystal quality.
A Deep Ensemble (DE) of Neural Networks serves as a surrogate model, trained on numerical data from a Finite Element Model (FEM). Optimization is carried out using NSGA-II and NSGA-III, two variants of Genetic Algorithms that explore trade-offs between competing objectives and identify high-performing candidate solutions. Results show that NSGA-II outperforms NSGA-III.
The optimal solutions correctly captured known trends, such as correlations between crystal size, pulling rate, and thermal stress.
A subset of the more intricate solutions was validated through new simulations, showing excellent prediction performance.
However, candidate solutions must still be verified by the FEM prior to experimental validation. 
This framework establishes a foundation for systematic, data-driven process optimization in FZ growth and can be extended to accelerate improvements in other crystal growth methods.
\end{abstract}

\section{Introduction}
\label{sec:intro}
Among the techniques for growing silicon (Si) crystals from the melt, \gls{FZ} stands out for its ability to produce ultra-pure, dislocation-free single crystals \cite{MUIZNIEKS2015}.
These crystals are crucial for high-power electronics in the green energy sector and advanced detection systems used in cutting-edge research, such as the Einstein Telescope \cite{Schnabel2010, Vieira2024}.
In inductively heated \gls{FZ}, a molten zone is sustained by a high-frequency induction coil as it moves vertically through a polycrystalline Si feed rod.
The melt is supported solely by the solidified Si below.
Surface tension stabilizes the free surface against \gls{EM} forces, gravity, hydrostatic pressure, and centrifugal force from rotation.
Crystal growth begins when a Si seed is brought into contact with the melt, establishing the orientation for the solidifying crystal lattice.
The thermal shock at this stage generates dislocations, which are line defects in the otherwise regular atomic arrangement.
These dislocations are subsequently eliminated using the Dash-neck technique \cite{Dash1959}.
As a crucible-free growth method, \gls{FZ} completely eliminates one major source of impurities, thereby allowing crystals reaching a purity of \SI{99.9999999}{\percent}.

\gls{FZ}-Si crystals up to \SI{200}{\milli\meter} in diameter can be grown industrially, and the \gls{IKZ} is developing its own process at this scale.
However, reliably attaining this size and scaling it up is notoriously challenging.
Large crystals increase the risk of process instabilities, including dislocation generation, electrical arcing, spiral growth, and melt spillage.
At the same time, crystals must meet stringent quality criteria for electronic applications, such as low point defect concentrations and uniform resistivity.
The complex interplay of physical phenomena in \gls{FZ} makes it extremely difficult to optimize all aspects of the process simultaneously.
Improving one aspect often compromises another, making this fundamentally a \gls{MOO} problem with productivity and crystal quality objectives.

Numerical simulations have long supported \gls{FZ} optimization \cite{Coriell1976, Coriell1977, DaRiva1978, Duranceau1986, Lie1989, Lie1990, Young1990, Muehlbauer1993, Muehlbauer1995, Ratnieks1999, Rudevičs2004, Dadzis2005, Wunshcer2014, Han2020, Surovovs2022, Ai2023, Vieira2024, Tsiapkinis2024}.
However, past studies optimized only a few parameters at a time.
Despite great advancements in computational power and numerical methods in the last decades, long computation times remain a hurdle.
Regardless of this limitation, \gls{MOO} inherently requires a substantial number of evaluations to effectively explore new parameter spaces that balance competing objectives.
To date, there appears to be no study that has conducted a systematic \gls{MOO} investigation for \gls{FZ} growth.

The growing adoption of data-driven and \gls{ML} methods presents new opportunities for optimizing crystal growth.
Various techniques have recently been applied to processes such as Czochralski, Vertical Gradient Freeze, and Directional Solidification \cite{Dropka2014, Dropka2017, Boucetta2019, Kutsukake2020, Dropka2021, Dropka2022, Kutsukake2022, Ghritli2022, Tang2023, Takehara2023, Chappa2025}.
In the case of \gls{FZ}, \gls{ML} has primarily been used to support process control and dynamics \cite{Omae2022, Tosa2023, Chen2024}.
The first exploratory study examining the interplay between multiple parameters and outcomes in \gls{FZ} was presented in \cite{Vieira2025}.

Data-driven techniques can effectively tackle \gls{MOO} problems as well.
One approach is the surrogate-based optimization \cite{Queipo2005}, where a surrogate model replaces costly experiments or simulations and is repeatedly queried by an optimization algorithm.
Commonly, \glspl{NN} serve as surrogates, learning patterns from past data to make fast, accurate predictions \cite{Uhrig1995}, while \glspl{GA} optimize complex problems using biologically inspired strategies \cite{Marler2004}.
In bulk crystal growth, this method was applied to the Top Seeded Solution Growth process \cite{Isono2022}.
By generating multiple optimal parameter combinations, this framework reveals inherent trade-offs between competing objectives.

To gain insights into the combined optimization of multiple process parameters in \gls{FZ} growth, specifically the geometry of the inductor and the growth conditions, we employ a surrogate-based framework that combines \glspl{NN} with \glspl{GA}.
We target eight representative objectives related to productivity and process stability.
Following the methodology in \cite{Vieira2025}, numerical simulations are performed across a wide range of parameter combinations to generate the necessary training data.
A surrogate model is then built using an ensemble of \glspl{NN}, referred to as a \gls{DE}, which enhances prediction robustness and accuracy \cite{Lakshminarayanan2017, Aggarwal2018_4}.
For the optimization task, we adopt the \gls{NSGA} family of \glspl{GA}, comparing the performance of NSGA-II \cite{Deb2002} and NSGA-III \cite{Deb2014} variants.
This approach yields a diverse set of candidate solutions that perform well across one or more objectives. By visualizing the trade-offs among objectives, the framework provides guidance for future numerical and experimental investigations.
To the best of the authors’ knowledge, this is the first time \gls{MOO} optimization was applied to \gls{FZ} crystal growth with more than three objectives, and where a \gls{DE} is used as a surrogate model in place of a single \gls{NN}.

\section{Methodology}

\subsection{Framework for Multi-Objective Optimization}
\label{sec:framkework}
\begin{figure}
    \centering
    \includegraphics[width=0.75\linewidth]{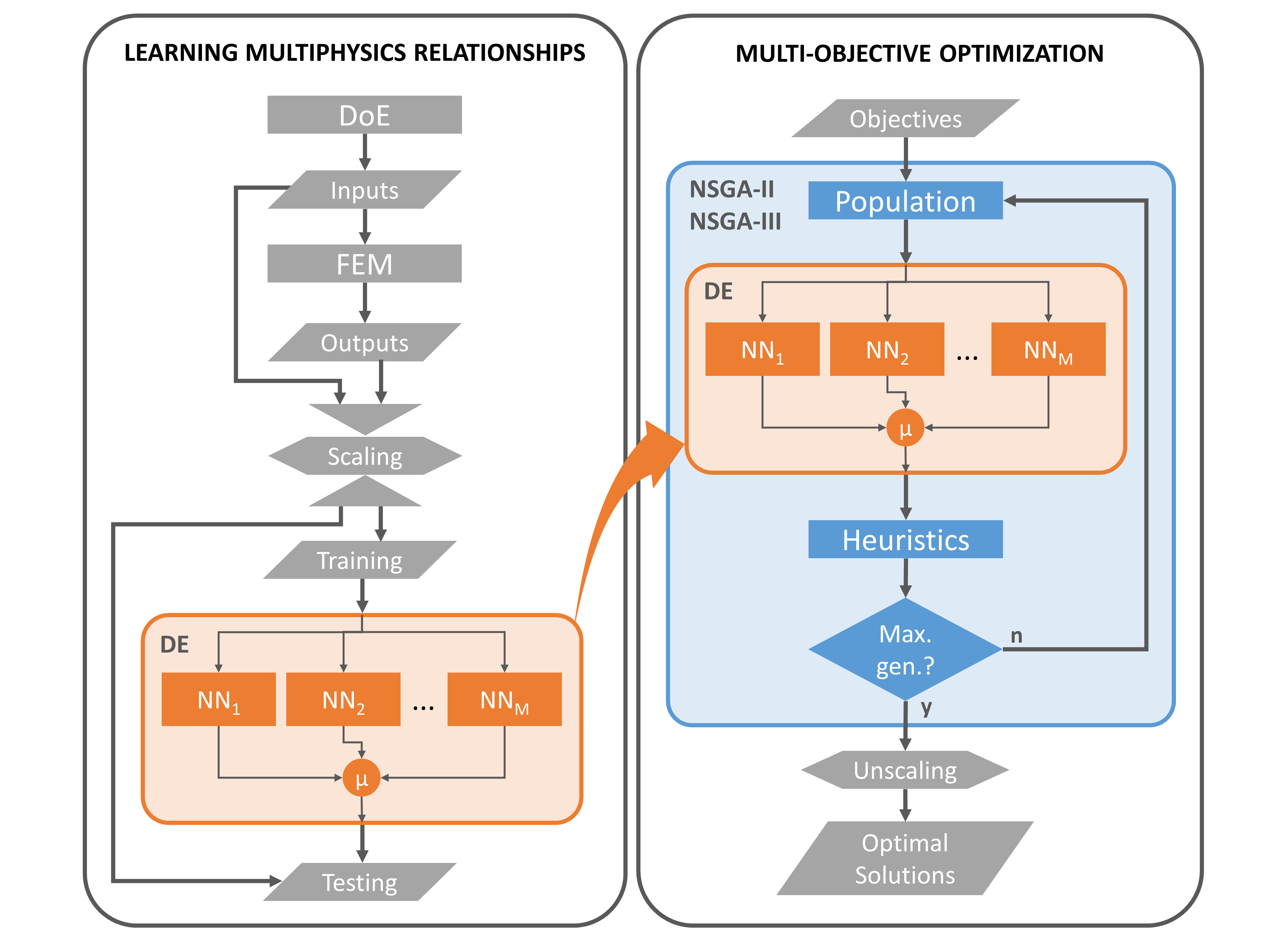}
    \caption{Framework for \gls{MOO} in \gls{FZ} growth.}
    \label{fig:MOOframework}
\end{figure}

The framework that reveals trade-offs between competing objectives in \gls{FZ}-Si growth is depicted in Figure~\ref{fig:MOOframework}. 
It is a combination of the following items:
The list of objectives (Section~\ref{sec:ObjectiveList});
A \gls{FEM} provides a physics-aware relationship between the model inputs and outputs by solving relevant governing equations of the \gls{FZ} process (Section~\ref{sec:numerical_model});
The \gls{DoE} selects input combinations such that the simulation results are informative (Section~\ref{sec:doe});
The \gls{DE} acts as a fast surrogate by learning to represent the functional dependence between the inputs and outputs from the \gls{FEM} (Section~\ref{sec:MLtechniques});
Finally, the \gls{GA} drives the \gls{DE} in seek of solutions that optimize the preset objectives (Section~\ref{sec:multiObjectiveOptimization}).

\subsection{List of Objectives}
\label{sec:ObjectiveList}
The eight objectives addressed in this work were selected based on their relevance to achieving a stable process and high crystal quality, as supported by both practical experience and literature:

\paragraph{Crystal Radius and Pulling Rate}
Maximizing either of these parameters improves productivity.
Therefore, understanding their interdependence is desirable.

\paragraph{Radial Thermal Gradient}
Spiral growth is the phenomenon where the crystal loses cylindrical structure and starts to grow in an unwanted helical shape \cite{Schwabe2011}.
The likelihood of spiral growth can be minimized by a sufficiently high radial thermal gradient $\radialThermalGradient$ at the \gls{TPL} (see TPL in Figure~\ref{fig:model}).

\paragraph{Interface Deflection}
The shape of the crystallization interface is notably affected by geometrical and process parameters such as the geometry of the inductor, reflector, pulling rate, \textit{etc}.

\paragraph{Exceed Stress}
High mechanical stress caused by thermal expansion gradients can spontaneously generate dislocations.
In theory, dislocations appear when the Von Mises stress $\stress_\text{vM}$ exceeds the temperature-dependent critical resolved shear stress $\stress_\text{CRSS}(T)$.
However, simulations of successful \gls{FZ} growth processes indicate that crystals tolerate stresses several tens of \SI{}{MPa} above this threshold.
To better quantify the risk of dislocation generation, Muižnieks \textit{et al.} introduced the exceed stress variable $\exceedStress$, defined as the difference between $\stress_\text{vM}$ and $\stress_\text{CRSS}(T)$ \cite{Muiznieks2001}.
Specifically, higher values of $\exceedStress$ correspond to increased risk of dislocation generation.

\paragraph{Voltage} 
Electrical arcing disrupts \gls{FZ} growth when the voltage drop across the inductor terminals $\voltage$ exceeds the breakdown voltage $\voltage_\text{b}$ of the surrounding gas.
Paschen's law models $\voltage_\text{b}$ taking into account gas composition \cite{Kuffel2000}.
While somewhat high $\voltage_\text{b}$ is possible depending on gas pressure and gap distance between terminals, a practical analysis requires experimentally determining material constants.
For sake of simplicity, only $\voltage$ is analyzed in this study.

\paragraph{EM Homogeneity}
Dislocations can also be induced if newly crystallized regions are remelted.
This may happen due to fluctuations in the induced \gls{EM} power $\inducedHeatEM$ below the main slit, in combination with crystal pulling rate and rotation \cite{Menzel2013}.
To evaluate the impact of the main slit, a helper variable $\EMhomogeneity$ was introduced:
\begin{equation}
    \EMhomogeneity = \oint_{\triplePhaseLine} |\inducedHeatEM - \langle \inducedHeatEM \rangle| \,ds \bigg/ \triplePhaseLine \text{.}
\end{equation}
The integral is evaluated along the closed \gls{TPL} path $\triplePhaseLine$, where $\langle \inducedHeatEM \rangle$ is the averaged $\inducedHeatEM$ along $\triplePhaseLine$.
Minimizing $\EMhomogeneity$ corresponds to increasing the homogeneity of $\inducedHeatEM$ along \gls{TPL}.

\paragraph{Voronkov's Criterion}
The presence of point defects in the Si lattice, e.g. vacancies and self-interstitials, are detrimental to electronic performance as they interfere with carrier lifetime.
The Voronkov's criterion $\VoronkovCriterion$ indicates if these defects are thermodynamically favorable \cite{Voronkov1982}:
\begin{equation}
    \VoronkovCriterion = \frac{\crystalPullingRate}{\thermalGradient} \text{,}
    \label{eq:VoronkovsRatio}
\end{equation}
with $\crystalPullingRate$ the crystallization rate and $\thermalGradient$ the axial thermal gradient.
At a critical value $\VoronkovCriterion = \VoronkovCriterion_\text{crit}$, both kinds of defects annihilate each other.
If $\VoronkovCriterion < \VoronkovCriterion_\text{crit}$, the crystal is self-interstitial rich, whereas if $\VoronkovCriterion > \VoronkovCriterion_\text{crit}$ it is vacancy-rich.
Empirically obtained values for $\VoronkovCriterion_\text{crit}$ in \gls{FZ}-Si range between $\VoronkovCriterion_\text{crit}^\text{LB} =$ \SI{1.3e-3}{} and $\VoronkovCriterion_\text{crit}^\text{UB} =$ \SI{2.2e-3}{\square\centi\meter \per\minute \per\kelvin} \cite{Vanhellemont2015}.

\subsection{Numerical Model}
\label{sec:numerical_model}
The utilized \gls{FEM}, implemented in \texttt{COMSOL Multi\-physics\textsuperscript{\textregistered{}}} \cite{comsolv61}, is detailed in \cite{Vieira2025} with Si and Cu material properties consistent with \cite{Vieira2024}.
The model consists of two coupled components, one in 3D for the \gls{EM} field and the other 2D-axisymmetric for heat transfer and stress field.
It assumes steady-state growth of Si crystals in the cylindrical phase.
Thus, despite crystal diameter $d$ not being a free parameter in the \gls{FZ} process, the \gls{FEM} considers crystal radius $r = d/2$ as an input.

Regarding the governing equations, the \gls{FEM} solves Maxwell equations in the gas domain, heat transfer within Si domains with radiation exchange between Si and inductor surfaces, and the von Mises stress within crystalline Si (see \cite{Wunshcer2014, Vieira2025} for more details).
Melt flow is neglected as it has a limited influence on the set of objectives investigated in this work (implications of this approximation can be found in \cite{Ratnieks2007}).
As in \cite{Vieira2025}, the external phase boundaries are fixed, and the boundaries between Si phases are computed with the Apparent Heat Capacity formulation.

\begin{figure}
    \centering
    \includegraphics[width=0.75\linewidth]{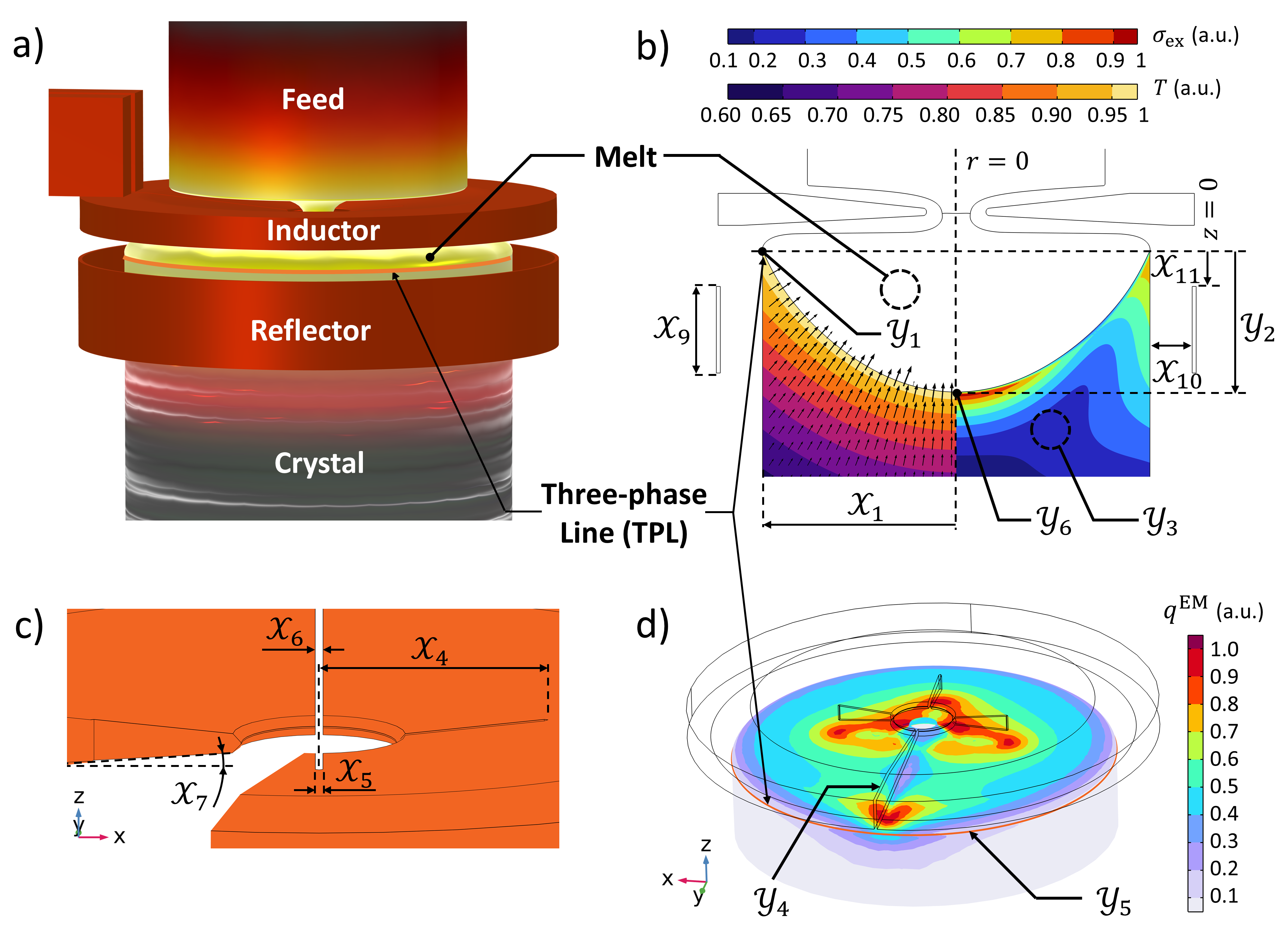}
    \caption{
    Depiction of the \gls{FEM}.
    a) Main constituents of the \gls{FZ} process.
    b) Inputs and outputs in the 2D-axisymmetric model: temperature and its gradient (left), stress field (right).
    c) Inductor geometrical parameters.
    d) \gls{EM} quantities in the 3D model.
    }
    \label{fig:model}
\end{figure}

The model was parametrized on the \num{12} inputs $\modelInput$ from Table~\ref{tab:modelInputs}, selected based on both expert knowledge and feature importance analysis from \cite{Vieira2024} for the outputs of interest.
Simulations were performed with input combinations within the specified ranges (Section~\ref{sec:doe}), from which six representative outputs $\modelOutput$, also shown in Table~\ref{tab:modelInputs}, were extracted.
All outputs and the geometric inputs are illustrated in Figure~\ref{fig:model}.

\begin{table}
    \centering
        \begin{tabular}{clll}
            \hline
             & \textbf{Description} & \textbf{Range} & \textbf{Unit} \\
             \hline
             \multicolumn{2}{l}{\textbf{Inputs} $\modelInput$}\\
            $\inputCrystalRadius$       & Crystal radius                                          & $\interval{75}{100}$                   & \SI{}{\milli\meter} \\
            $\inputPullingRate$         & Crystal pulling rate                                    & $\interval{1.0}{3.5}$                  & \SI{}{\milli\meter \per\minute} \\
            $\inputSideSlits$           & Inductor side slits (quantity of)                       & \{\num{2}, \num{3}, \num{4}, \num{5}\} & \SI{1}{} \\
            $\inputSideSlitLength$      & Inductor side slits length                              & $\interval{30}{60}$                    & \SI{}{\milli\meter} \\
            $\inputSideSlitWidth$       & Inductor side slits width                               & $\interval{1.5}{4.0}$                  & \SI{}{\milli\meter} \\
            $\inputMainSlitWidth$       & Inductor main slit width                                & $\interval{0.5}{4.0}$                  & \SI{}{\milli\meter} \\
            $\inputBottomAngle$         & Inductor bottom angle                                   & $\interval{1.0}{4.0}$                  & \SI{}{\degree} \\
            $\inputFrequency$           & Inductor frequency                                      & $\interval{2.0}{3.5}$                  & \SI{}{\mega\hertz} \\
            $\inputReflectorHeight$     & Reflector height                                        & $\interval{20}{80}$                    & \SI{}{\milli\meter} \\
            $\inputReflectorRadius$     & Reflector radius (offset from $\inputCrystalRadius$)    & $\interval{5}{40}$                     & \SI{}{\milli\meter} \\
            $\inputReflectorPosition$   & Reflector position (relative to \gls{TPL})                    & $\interval{-5}{10}$                    & \SI{}{\milli\meter} \\
            $\inputReflectorEmissivity$ & Reflector emissivity                                    & $\interval{0.10}{0.85}$                & \SI{1}{} \\
            \hline
            \multicolumn{2}{l}{\textbf{Outputs} $\modelOutput$}\\
            $\outputRadialThermalGradient$    & Radial Thermal Gradient at \gls{TPL}             & -                                         & \SI{}{\kelvin\per\centi\meter} \\
            $\outputDeflection$               & Deflection of the crystallization interface             & -                                         & \SI{}{\milli\meter} \\
            $\outputExceedStress$             & Exceed Stress in the crystal                            & -                                         & \SI{}{\mega\pascal} \\
            $\outputVoltage$                  & Voltage drop across inductor terminals                  & -                                         & \SI{}{\volt} \\
            $\outputEMhomogeneity$            & Homogeneity of the induced \gls{EM} power at \gls{TPL}        & -                                         & \SI{}{\square\watt \per\meter} \\
            $\outputVoronkovCriterion$        & Voronkov's Criterion at the center of the crystal       & -                                         & \SI{}{\square\centi\meter \per\minute \per\kelvin} \\
            \hline
        \end{tabular}
    \caption{Twelve selected model inputs with their ranges, and the six representative outputs.
    All parameters are continuum values except for $\inputSideSlits$, which can only be a positive integer.}
    \label{tab:modelInputs}
\end{table}

In order to provide point-wise training data for the \glspl{NN} in the subsequent step, we introduce the following definitions:
Due to the Apparent Heat Capacity formulation, there is no sharp boundary between the Si phases in our \gls{FEM} \cite{Vieira2025}.
Thus, $\outputVoronkovCriterion$ is the averaged $\VoronkovCriterion$ at $r=0$ in the crystallization region.
$\outputDeflection$ is the $z$ coordinate of the crystallization interface $\deflection$, where $T$ equals melting temperature of Si.
$\outputExceedStress$ is the maximum value of $\exceedStress$ within the crystal.
$\outputRadialThermalGradient$ is the averaged Radial Thermal Gradient in a circular region with a radius of \SI{0.5}{\milli\meter} and centered at \gls{TPL}.
Finally, $\outputEMhomogeneity = \EMhomogeneity$ and $\outputVoltage$ is a global variable in the \gls{FEM}.

\subsection{Design of Experiments}
\label{sec:doe}
In bulk crystal growth, each output is affected by a complex combination of inputs \cite{Vieira2025}.
\gls{DoE} techniques allow to consider these interactions while requiring as few experiments as possible.
In light of this, we use \gls{LHS} to sample input combinations for our simulations \cite{Mckay1979}.

\gls{LHS} stratifies each dimension from a space of $\dimensionsLHS$ inputs into $\samplesLHS$ intervals, creating $\samplesLHS^{\dimensionsLHS}$ subsets.
Out of these, $\samplesLHS$ random points, each representing an input combination, are sampled without subset repetition.
This ensures that each input dimension is sampled uniformly while preserving randomness.
An appropriate coverage the $\dimensionsLHS$-dimensional input space depends on the number of samples, i.e. the more points, the less likely interaction effects will be missed.
For our study, $\samplesLHS =$ \num{2500} simulations with input combinations within the ranges from Table~\ref{tab:modelInputs}.

\subsection{Machine Learning Techniques}
\label{sec:MLtechniques}

\subsubsection{Artificial Neural Networks}
The strong nonlinear multiphysics of \gls{FZ} growth was modeled with \glspl{NN}, which can approximate any continuous function \cite{Uhrig1995}.
Given a dataset containing $\samplesLHS$ samples of a multivariate relationship between a vector of inputs and a vector of outputs, $\trainingData = \{\mathbf{x}, \mathbf{y}(\mathbf{x})\}_{i=1}^S$, a feed-forward \gls{NN} extracts hierarchical representations from the training data.
The approximation $\hat{\mathbf{y}}(\mathbf{x}) \approx \mathbf{y}(\mathbf{x})$ is obtained by computing $L$ compositions of transformations as follows:
\begin{equation}
    \hat{\mathbf{y}} = \left( \mathbf{a}^{[L]} \circ \mathbf{a}^{[L-1]} \circ \dots \circ \mathbf{a}^{[1]} \right) (\mathbf{x}) \text{,}
    \label{eq:NN}
\end{equation}
where $\mathbf{a}^{[l]}$ is the vector of nonlinear transformations of the $l$-th layer, also called neurons:
\begin{equation}
    \mathbf{a}^{[l]} = \pmb{\phi}^{[l]} \left( \mathbf{W}^{[l]}\mathbf{a}^{[l-1]} + \mathbf{b}^{[l]} \right) \text{.}
    \label{eq:NNneurons}
\end{equation}
Here, $\mathbf{W}^{[l]}$ and $\mathbf{b}^{[l]}$ denote the weight matrices and bias vectors (linear part), while $\pmb{\phi}^{[l]}(\cdot)$ is a vector of activation functions (nonlinear part).
In case of regression tasks, the last layer $\pmb{\phi}^{[L]}(\cdot)$ typically performs linear transformations to ensure continuous outputs.

The number of hidden layers ($L-1$), the number of neurons in each layer $\mathbf{a}^{[l]}$ and their respective activation functions $\pmb{\phi}^{[l]}(\cdot)$, are determined prior to training.
The $\{\mathbf{W}^{[l]}, \mathbf{b}^{[l]}\}$ parameters are then automatically adjusted during training to minimize some kind of error metric, called loss function $\loss$.
Modern \gls{ML} practice relies on gradient-based optimization, where backpropagation is used to compute $\loss$ w.r.t. $\{\mathbf{W}^{[l]}, \mathbf{b}^{[l]}\}$ gradients.
The updates are iteratively performed over $N$ passes through the training data, with a learning rate $\learningRate$ that determines the step size of each adjustment.

In our study, we implemented \glspl{NN} on \texttt{Tensorflow} \cite{tensorflow2016} with $\phi = \text{\gls{ReLU}} = \text{max}(0,x)$, Adam optimizer and $\learningRate = 0.001$ for parameter updates, $N =$ \num{100} and \gls{MSE} as $\loss$. The loss of an output $j$ is:
\begin{equation}
   \loss_j = \text{MSE}_j = \frac{1}{\samplesLHS} \sum_{i=1}^{\samplesLHS} \left( \hat{\mathbf{y}}_j^{(i)} - {\mathbf{y}}_j^{(i)}\right)^2 \text{,}
   \label{eq:loss_mse}
\end{equation}
and the overall \gls{NN} loss for our six outputs is $\loss = \sum_{j=1}^6 \loss_j$.
In order for $\loss$ to be meaningful, the columns of $\trainingData$ must share the same relative scale.

\subsubsection{Hyperparameter Optimization}
To determine the optimal architecture (number of hidden layers and neurons per layer), a hyperparameter optimization strategy using $k$-fold cross-validation was employed.
Specifically, $\trainingData$ was divided into \SI{90}{\percent}-\SI{10}{\percent} for training and testing, $\trainingData_\text{train}$ and $\trainingData_\text{test}$ respectively, where each column was scaled between $\interval{0}{1}$.
$\trainingData_\text{train}$ was further divided into \num{10} folds, where \num{9} were used for training and \num{1} for preliminary validation.
Then, \num{1000} architectures with 1-10 hidden layers and 2-64 neurons were tested on each $\trainingData_\text{train}^k$ fold \num{10} times, each time using a different fold for validation.
The search for architectures was driven by the \gls{TPE} sampling algorithm \cite{Watanabe2023}.
In the end of the search, the architecture with the lowest average $\loss$ was retrained on the full $\trainingData_\text{train}$, this time using the unseen $\trainingData_\text{test}$ to access prediction performance.

\subsubsection{Deep Ensembles}
Once trained, all $\{\mathbf{W}^{[l]}, \mathbf{b}^{[l]}\}$ parameters are fixed and the \gls{NN} can be used to infer outcomes based on new input sets.
However, since training a \gls{NN} is an inherently stochastic process due to random initialization of these parameters \cite{Aggarwal2023_1}, slightly different predictions would be obtained if the \gls{NN} was retrained under identical circumstances.
Regardless of how small the difference, this stochastic variability may impact reliability in our setting, as most relationships in \gls{FZ} are subtle \cite{Vieira2025}.
To obtain more consistent predictions, we use a \gls{DE}, which is an Ensemble Learning technique for \glspl{NN} \cite{Lakshminarayanan2017, Aggarwal2018_4}.
The idea is to train $M$ models $(\hat{\mathbf{y}})_m $ with the same architecture, but with random initialized $\{\mathbf{W}^{[l]}, \mathbf{b}^{[l]}\}$, and average their predictions as a mixture:
\begin{equation}
    \text{\gls{DE}} = \frac{1}{M} \sum_{m=1}^{M} \left( \hat{\mathbf{y}} \right)_m \text{.}
\end{equation}
Other than that the \gls{DE} behaves exactly as a single \gls{NN}.

In addition to \gls{MSE}, the coefficient of determination for each output $\coeffDetermination_j$ was also reported:
\begin{equation}
   \text{R}^2_j = \frac{\sum_{i=1}^{\samplesLHS} \left( \text{\gls{DE}}_j^{(i)} - \bar{\mathbf{y}}_j \right)^2}{\sum_{i=1}^{\samplesLHS} \left( \mathbf{y}_j^{(i)} - \bar{\mathbf{y}}_j \right)^2} \text{,}
   \label{eq:r2}
\end{equation}
with $\bar{\mathbf{y}}_j = \sum_{i=1}^{\samplesLHS} \mathbf{y}_j^{(i)}$.

\subsection{Multi-Objective Optimization}
\label{sec:multiObjectiveOptimization}
\subsubsection{General Definitions}
Solving \gls{MOO} problems translates into finding one or multiple solutions that balance competing objectives in a set $\objective$, possibly subjected to constraints \cite{Coello2018, Marler2004}.
A solution $i$ consists of an input vector $\mathbf{x}^{(i)} = \{x_1^{(i)}, x_2^{(i)}, \dots\}$ which maps to an output vector $\mathbf{y}^{(i)} = \{y_1^{(i)}, y_2^{(i)}, \dots\}$, where each element in $\mathbf{y}^{(i)}$ is a function of all elements from $\mathbf{x}^{(i)}$.
In general, objectives can be expressed as functions and penalty terms depending on both inputs and outputs, i.e. $\objective_1 = f_1(\mathbf{x}, \mathbf{y}) + \epsilon_1(\mathbf{x}, \mathbf{y})$.
A prototypical \gls{MOO} problem takes the form:
\begin{equation}
    \begin{aligned}
    \text{Minimize } & \mathbf{\objective} = \{ \objective_1, \objective_2, \dots, \objective_o \}\\
    \text{such that} &
    \begin{cases}
        I_p(\mathbf{x}, \mathbf{y}) \le 0 \text{, } &p = \text{1, 2,} \dots, q\\
        E_r(\mathbf{x}, \mathbf{y}) = 0 \text{, } &r = \text{1, 2,} \dots, s\\    
    \end{cases}
    \end{aligned}
\label{eq:MOO_equation}
\end{equation}
with $o \ge 2$ objectives, $q$ inequality constraints $I_p$, and $s$ equality constraints $E_r$.
Constraints and penalties are case-specific tailored to guide solutions during optimization.
A given objective can be maximized by minimizing its negative value.

There are many approaches to frame and solve Equation~\ref{eq:MOO_equation}, and a distinction between \textbf{a priori} and \textbf{a posteriori} methods must be made \cite{Marler2004}.
\textbf{A priori} methods require further information about the objectives before optimization starts.
For example, a \gls{MOO} problem can be transformed into a single global optimization problem by constructing a weighted sum with all objectives (scalarization method), or the user can rank the objectives to be optimized in order of importance (Lexographic method) \cite{Marler2004}.
Such methods are useful when there are relatively few objectives with known trade-offs and meaningful interpretations for the weights and ranking.
In a more general case for exploratory assessment of many trade-offs, \textbf{a posteriori} methods are preferred because they make no assumptions and each objective is optimized individually.
The scope of our work falls into the second case and the algorithms used are explained in Section~\ref{sec:genetic_algorithm}.

\subsubsection{Pareto Optimal Solutions}
\label{sec:paretoOptimal}
After optimization, the performance of solutions can be assessed in a plot.
Let us illustrate this process by considering $o=4$, where we seek to maximize \{$\objective_1 = x_1$, $\objective_2 = y_1$\}, and minimize \{$\objective_3 = x_2$, $\objective_4 = y_2$\}.
It should be noted that, although optimizing inputs is trivial, i.e. just maximize or minimize their values, the effect on the outputs must be analyzed.

\begin{figure}
    \centering
    \includegraphics[width=0.85\linewidth]{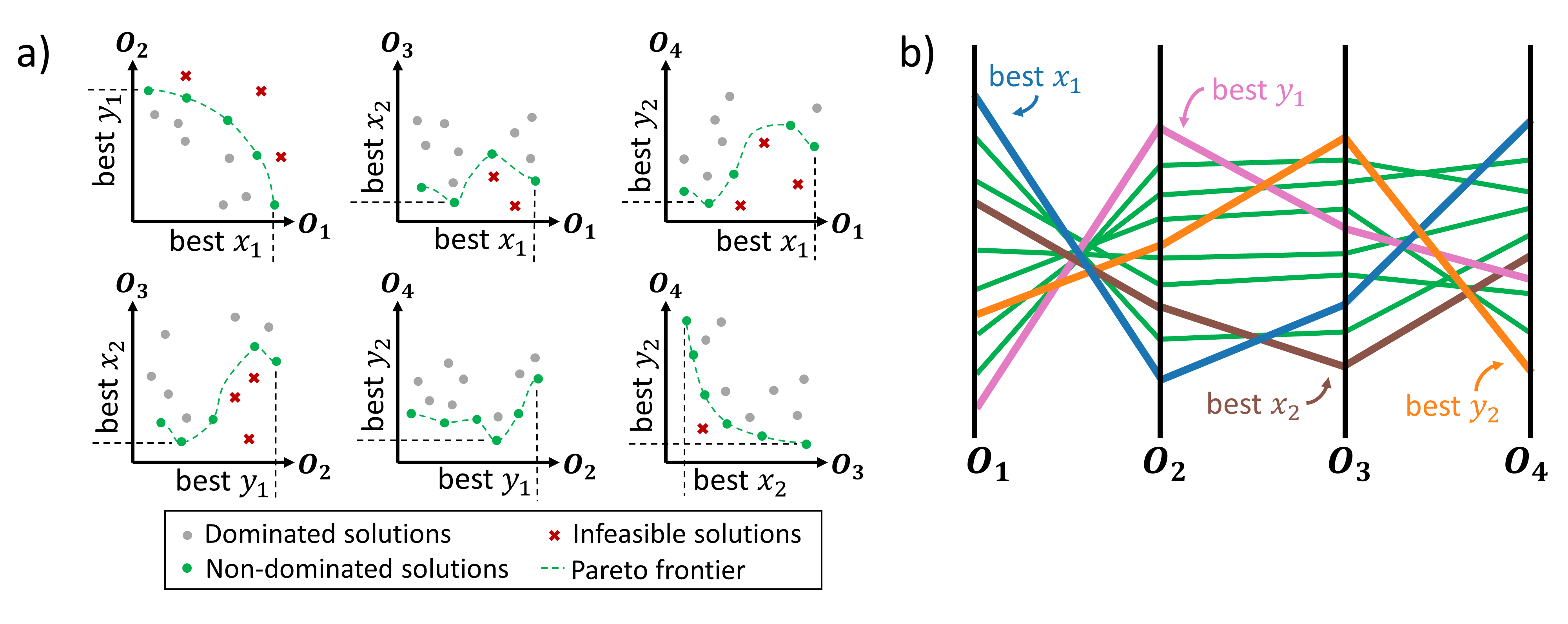}
    \caption{Optimal solutions. a) Two-dimensional scatter plots with all objective combinations. b) Parallel-coordinate plot with all objectives.}
    \label{fig:paretoExample}
\end{figure}

Figure~\ref{fig:paretoExample} shows two ways of visualizing optimal solutions.
In Figure~\ref{fig:paretoExample}a), scatter plots show all pairwise combinations of objectives.
A solution is considered non-dominated (or Pareto optimal) if no other solution improves all objectives simultaneously.
Unfeasible solutions violate one or more constraints, e.g. cannot be realized in practice.
The Pareto frontier separates dominated solutions from infeasible solutions.

The fact that most solutions appear dominated in Figure~\ref{fig:paretoExample}a) is misleading.
Since scatter plots are limited to three dimensions at most, information about any additional objective will be obscured.
Indeed, a solution belonging to a Pareto front according to a pair of objectives may perform poorly in another.
To overcome this limitation, Figure~\ref{fig:paretoExample}b) uses a parallel coordinate plot to show all objectives simultaneously.
Each line represents a feasible solution, revealing how it performs across every objective.
As a result, more solutions appear non-dominated when considered in the full objective space.

Trade-offs between multiple objectives become clearer using parallel coordinate plot:
\begin{itemize}
    \item A high number of crossing lines in X-shape between two axis indicate strong negative correlation of the objectives;
    \item Many parallel lines indicate strong positive correlation;
    \item A mixed pattern of the two above indicate a weak correlation.
\end{itemize}
Thus, in Figure~\ref{fig:paretoExample}b), solutions with high $\objective_1$ have in general low $\objective_2$ and $\objective_3$, and might have high or low $\objective_4$.
Understanding these trade-offs helps the user narrow down a suitable subset of solutions.

\subsubsection{Genetic Algorithms}
\label{sec:genetic_algorithm}
To obtain Pareto-optimal solutions discussed in Section~\ref{sec:paretoOptimal}, \glspl{GA} were used.
\glspl{GA} are a family of algorithms that apply heuristics inspired in biological evolution to solve optimization problems \cite{Marler2004, Kramer2017, Coello2018}.
In \gls{GA} terminology, a given solution is referred to as an \textit{individual}, and many individuals form a \textit{population} of candidate solutions.
Each individual is evaluated according to how well it performs across the set of objectives $\objective$.
By applying the evolutionary principles of \textit{selection}, \textit{crossover} and \textit{mutation} over many \textit{generations}, the \gls{GA} explores the space of available solutions in seek of optimal individuals.
Namely, the characteristics of the best evaluated individuals from the current generation are combined to produce new individuals in the next generation, thereby creating a natural-selection pressure.
After a sufficient number of generations or a convergence is reached, the surviving individuals will be considered Pareto-optimal.

To solve \gls{MOO} problems, standard \glspl{GA} must be adapted to not only favor Pareto-optimal solutions, but also ensure diversity, otherwise solutions may cluster around local or suboptimal region of the Pareto front.
Two widely used adaptations are \gls{NSGA}-II and \gls{NSGA}-III.
Both rank solutions by non-dominance, but differ in how they promote diversity:
\begin{itemize}
    \item \textbf{NSGA-II}:
    uses the Crowding Distance metric, which estimates how isolated a solution is by computing the average normalized distance to neighboring solutions along each objective axis \cite{Deb2002, Coello2018}.
    \item \textbf{NSGA-III}:
    introduces a set of uniformly distributed reference points in the objective space \cite{Deb2014, Coello2018}.
    During selection, solutions are associated with the nearest reference point, and the number of associated solutions is limited to maintain an even spread.
    The number of reference points $P$ grows combinatorially with the number of objectives $o$ and a user-defined granularity parameter $g$, as given by the binomial coefficient:
    \begin{equation}
    P = \binom{o + g - 1}{g} = \frac{(o + g - 1)!}{g!(o - 1)!} \text{.}
    \label{eq:refPoints}
    \end{equation}
\end{itemize}

The literature distinguishes between multi-objective ($1 < o < 4$) and many-objective ($o \ge 4$) problems, recommending \gls{NSGA}-II for the former and \gls{NSGA}-III for the latter \cite{Deb2014}.
This distinction arises because Crowding Distance loses efficiency in high dimensions.
Paradoxically, \gls{NSGA}-II has been found to outperform \gls{NSGA}-III in many-objective problems when objectives are correlated \cite{Hisao2016}.
This situation underscores the need for a comparative study in our setting.

The chosen objectives for \gls{FZ} are listed in Table~\ref{tab:objectives} and were evaluated with \gls{NSGA}-II and \gls{NSGA}-III.
The constraints in $\objective_1$ and $\objective_2$ define performance targets for a high-productivity regime.
A lower bound for $\objective_3$ is not formally established.
In $\objective_4$, $\objective_5$ and $\objective_6$, the constraints are informed estimates based on domain knowledge from past experiments.
The constraint in $\objective_7$ ensures physical results.
Finally, $\objective_8$ has a penalty term because constraining $\outputVoronkovCriterion$ within [$\VoronkovCriterion_\text{crit}^{LB}$, $\VoronkovCriterion_\text{crit}^{UB}$] would impose an overly restrictive condition.

As for the hyperparameters, we defined a population size of \SI{500}{} individuals, reproduction rate of \SI{70}{\percent}, mutation rate of \SI{5}{\percent}, and 250 generations.
In \gls{NSGA}-III, $g = 12$.
The code was implemented on \texttt{DEAP} library \cite{deap2012}.

\begin{table}
    \centering
        \begin{tabular}{lllrclc}
            \hline
             & \textbf{Objective} & \textbf{Quantity} & \multicolumn{3}{l}{\textbf{Constraint}} & \textbf{Penalty}\\
            \hline
            $\objective_1$ & Maximize & Crystal Radius & $\inputCrystalRadius$ & $\ge$ &  \SI{85}{\milli\meter} & - \\
            $\objective_2$ & Maximize & Pulling Rate & $\inputPullingRate $ & $\ge$ &  \SI{1.5}{\milli\meter \per\minute} & - \\
            $\objective_3$ & Maximize & Radial Thermal Grad. & $\outputRadialThermalGradient$ & &  (none)  & - \\
            $\objective_4$ & Minimize & Interface Deflection & $\outputDeflection $ & $\le$ &  \SI{70}{\milli\meter} & - \\
            $\objective_5$ & Minimize & Exceed Stress & $\outputExceedStress $ & $\le$ &  \SI{80}{\mega\pascal} & - \\
            $\objective_6$ & Minimize & Voltage &  $\outputVoltage $ & $\le$ &  \SI{1200}{\volt} & -\\
            $\objective_7$ & Minimize & EM Homogeneity & $\outputEMhomogeneity $ & $\ge$ &  \SI{0}{\watt \per\squared\meter} & -\\
            $\objective_8$ & Interval & Voronkov's Criterion & $\outputVoronkovCriterion$ & & (none) & 
                $\begin{cases}
                    (\VoronkovCriterion_\text{crit}^\text{LB} - \outputVoronkovCriterion)^*, & \text{if } \outputVoronkovCriterion < \VoronkovCriterion_\text{crit}^\text{LB} \\ 
                    (\outputVoronkovCriterion - \VoronkovCriterion_\text{crit}^\text{UB})^*, & \text{if } \outputVoronkovCriterion > \VoronkovCriterion_\text{crit}^\text{UB} \\ 
                    0, & \text{otherwise}
                \end{cases}$ \\
                \hline
        \end{tabular}
    \caption{Selected objectives for \gls{MOO} in \gls{FZ}.
    $(\cdot)^*$ indicates that the unit is disregarded.}
    \label{tab:objectives}
\end{table}

\section{Results and Discussion}
\subsection{Simulation Results}
Each simulation produces outcomes resembling the fields of Figure~\ref{fig:model}b and Figure~\ref{fig:model}d.
All \num{2500} simulations were visually inspected, and after removing the ones presenting nonphysical results, e.g. freezing of the melt core, \num{2390} simulations remained.
The distributions of inputs and outputs from the \gls{FEM} can be seen in Figure~\ref{fig:violinplots} (Section~\ref{sec:GAsolutions}).
The \gls{FEM} inputs exhibit uniform distributions, a characteristic of the Latin Hypercube Sampling method.
As for the \gls{FEM} outputs, $\modelOutput_1$ to $\modelOutput_5$ each display a single mode with some level of skewness, while $\modelOutput_6$ has three distinct modes.

\subsection{Deep Ensemble Surrogate Model}
\begin{figure}
    \centering
    \includegraphics[width=0.65\linewidth]{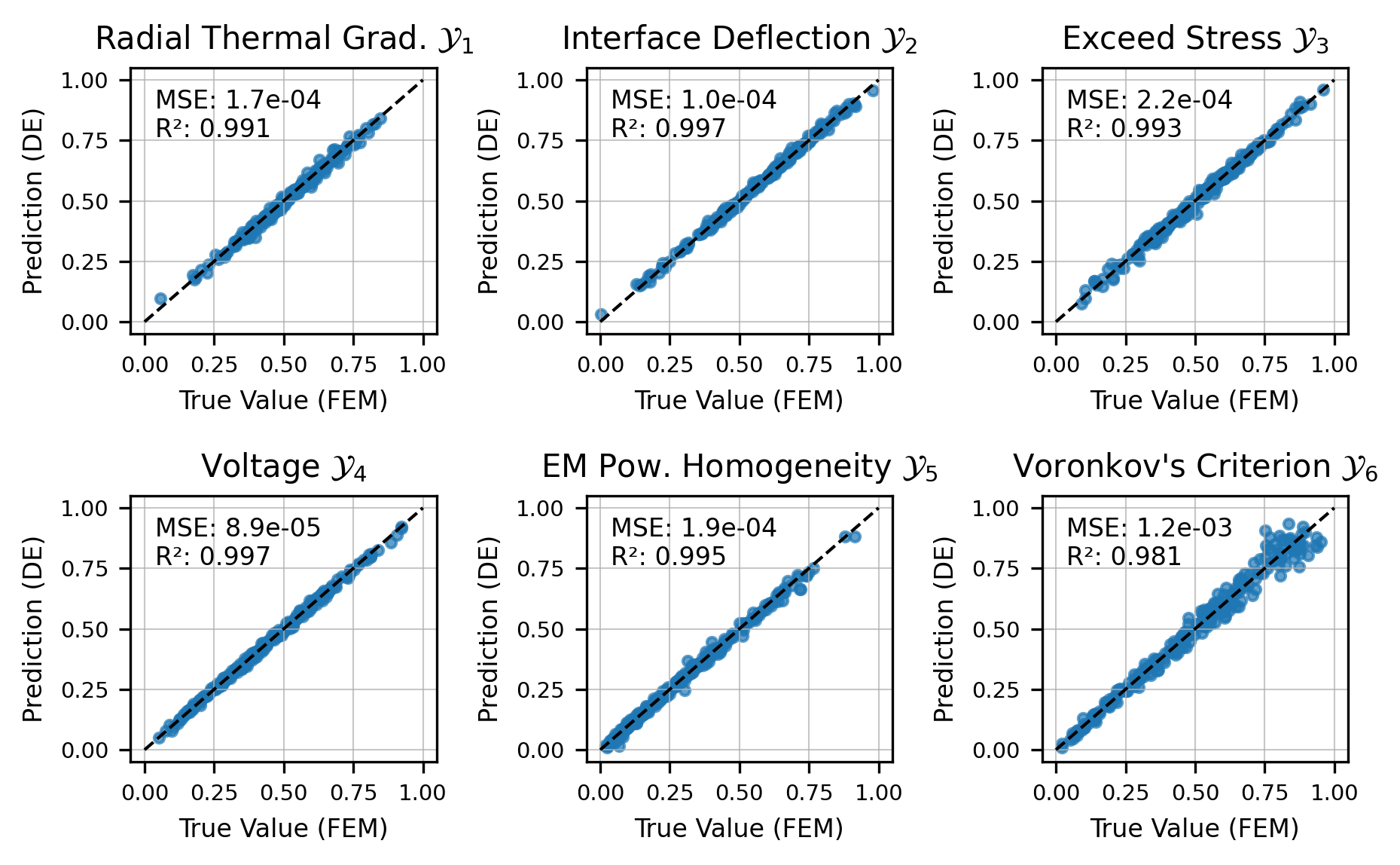}
    \caption{Parity plots and metrics, evaluated on $\trainingData_\text{test}$, of \gls{DE} predictions vs. the true \gls{FEM} simulations for each output.}
    \label{fig:parity_plots}
\end{figure}

The best architectures generated by the \gls{TPE} algorithm and their \glspl{MSE} for each count of hidden layers ($1,2,\dots,10$) are presented in Table~\ref{tab:bestNNs} (Appendix A).
To obtain the \gls{DE}, the \gls{NN} architecture with the lowest \gls{MSE} was retrained $M = 10$ times on the full $\trainingData_\text{train}$, this time using the unseen $\trainingData_\text{test}$ for final validation.
Figure~\ref{fig:parity_plots} shows parity plots for each output of the obtained \gls{DE} model, evaluated on the scaled $\trainingData_\text{test}$.
Visually, the dots cluster in the vicinity of the diagonal \gls{DE} = \gls{FEM}.
This is quantified by both $\text{\gls{MSE}}_j$ $\le \SI{1.2e-3}{}$ and $\coeffDetermination_j \ge 0.981$ for all outputs $j$.
The \gls{DE} functions as an 'advisory body' of specialists where each \gls{NN}, having learned a slightly different representation of the true relationship, will miss its target outputs by different amounts in different directions.
By averaging multiple predictions, the ensemble improves generalization and mitigates the risk of overfitting in the training data.

Noticeably in Figure~\ref{fig:parity_plots}, as $\outputVoronkovCriterion$ increases, the predictions deviate further from the \gls{DE} = \gls{FEM} line.
To analyze this behavior, it is reminded that $\VoronkovCriterion$ is defined as the ratio of crystallization rate to axial thermal gradient (Equation~\ref{eq:VoronkovsRatio}).
Since our model considers steady growth, the crystallization rate at $r=0$ equals pulling rate, $\crystalPullingRate = \inputPullingRate$.
Consequently, the dots with high values of $\outputVoronkovCriterion$ typically correspond to simulations with high values of $\inputPullingRate$.
It was observed in \cite{Vieira2025} that high values of $\inputPullingRate$ correspond to a higher variability of $\outputVoronkovCriterion$, following a similar scatter pattern to that of Figure~\ref{fig:parity_plots}.
The fact that \glspl{NN} could not qualitatively improve on this behavior indicates that it is probably a numerical and not a physical phenomenon, as \glspl{NN} are not supposed to fit pure noisy data.
Nevertheless, the \gls{DE} does provide a sufficiently good predictive model for our practical purposes and can be used as a surrogate to bypass the \gls{FEM}.

\begin{figure}
    \centering
    \includegraphics[width=0.85\linewidth]{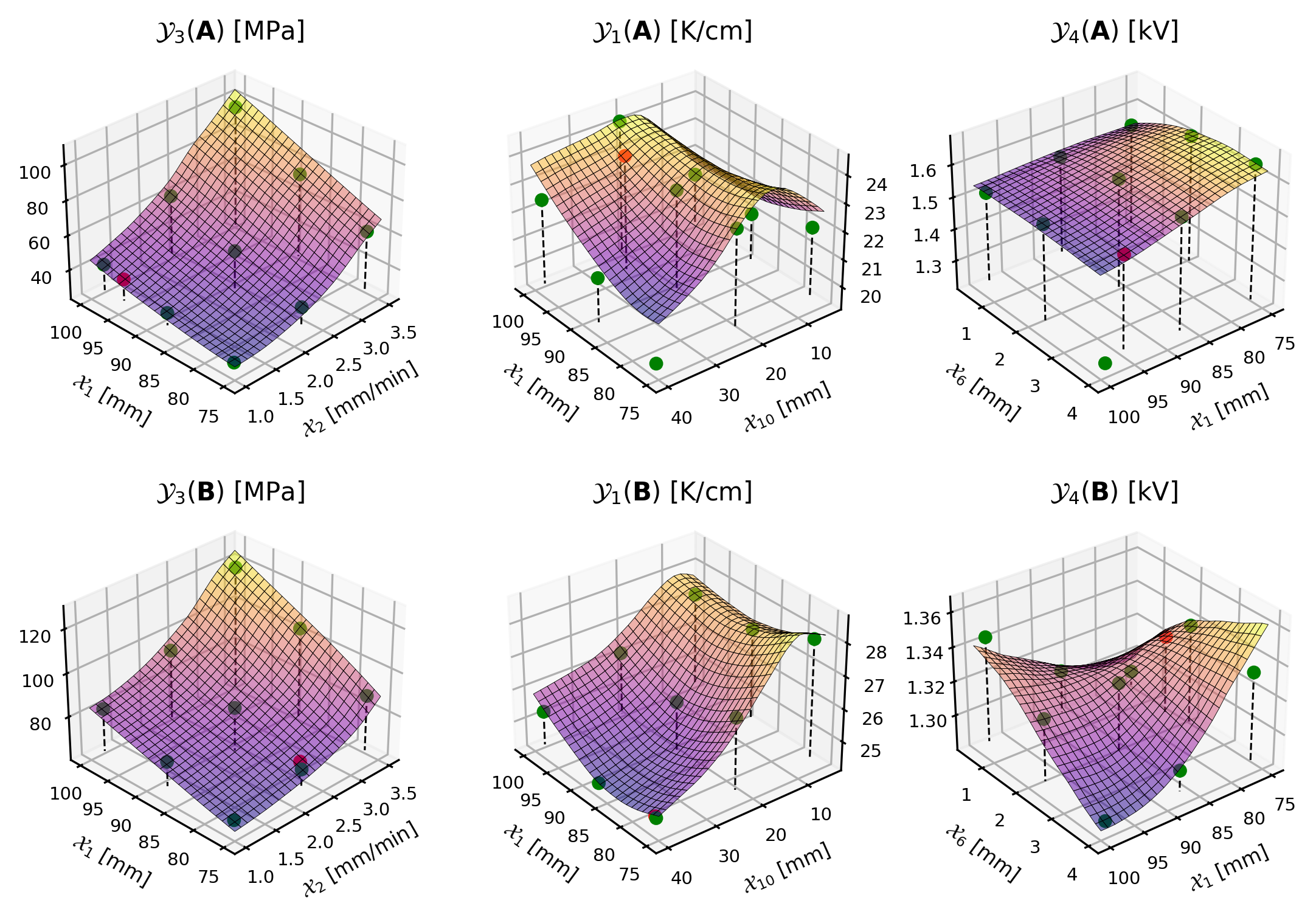}
    \caption{\gls{DE} surfaces: selected outputs as a function of two varying inputs with all other inputs fixed in combinations \textbf{A} and \textbf{B} (see Table~\ref{tab:combinationsSurfacePlot}).
    \gls{FEM} scatter plot: training data (red dots), and newly simulated data (green dots).}
    \label{fig:NN_surfacePlots}
\end{figure}

Being able to accurately and near-instantaneously predict simulation results enables us to shed light on complex relationships in \gls{FZ} growth.
To illustrate this, Figure~\ref{fig:NN_surfacePlots} shows selected outputs as a function of two varying inputs at a time, with all other inputs fixed in different combinations, denoted as cases \textbf{A} and \textbf{B}.
These combinations are detailed in Table~\ref{tab:combinationsSurfacePlot} found in Appendix B.

Among the illustrated outputs, only $\outputExceedStress$ exhibits a consistent functional dependence on $\inputCrystalRadius$ and $\inputPullingRate$ across both cases, with the difference between cases limited to a shift in magnitude.
In contrast, both $\outputRadialThermalGradient$ and $\outputVoltage$ display significant changes not only in magnitude but also in their functional relationships with the inputs, as evidenced by the distinct surface shapes in each case.
These variations are a manifestation of the intertwined relationships in \gls{FZ} growth.
Due to the sheer number of required queries to explore the space of possible input combinations, this task was offloaded to the \glspl{GA} (Section~\ref{sec:optimumSolutions}), which are capable of finding subspaces of interest more efficiently.

\subsubsection{Validation of the Deep Ensemble}
The red dots in Figure~\ref{fig:NN_surfacePlots} belong to the training data $\trainingData$.
The green dots were obtained by re-running the \gls{FEM} to compute $\modelOutput_i(\modelInput_j^\text{new}, \modelInput_k^\text{new}, \modelInput_{\cdots}^\trainingData)$, where $\modelInput_j^\text{new}$ and $\modelInput_k^\text{new}$ span a 3$\times$3 grid of the respective input ranges from Table~\ref{tab:modelInputs}.
All remaining input parameters $\modelInput_{\cdots}^\trainingData$ were fixed as one of the two instances from cases \textbf{A} or \textbf{B} (Table~\ref{tab:combinationsSurfacePlot} in Appendix B).
This selection provides insight into the true shape of the \gls{FEM} surface.

Despite slightly overestimating the true values, the \gls{DE} surfaces generally follow the shapes suggested by the dots remarkably well.
The obvious exception is the point ($\inputCrystalRadius = \SI{98.75}{\milli\meter}$, $\inputMainSlitWidth = \SI{3.825}{\milli\meter}$) in $\outputVoltage(\mathbf{A})$.
Incidentally, this is the closest point to the training case \textbf{A}.
While \glspl{NN} might trade off a poor performance on a few training examples for a better generalization, $\hat{\modelOutput}_4$ was the best performing output as from Figure~\ref{fig:parity_plots}.
The discrepancy is thus attributed to either an interaction effect not being properly covered by the \gls{DoE}, or a numerical artifact from the \gls{FEM}.
In $\outputVoltage(\mathbf{B})$, the greatest difference is \SI{20.2}{\volt} at ($\inputCrystalRadius = \SI{76.25}{\milli\meter}$, $\inputMainSlitWidth = \SI{3.89}{\milli\meter}$).

As for $\outputRadialThermalGradient(\mathbf{A})$, the green dot ($\inputCrystalRadius = \SI{98.75}{\milli\meter}$, $\inputReflectorRadius = \SI{22.5}{\milli\meter}$) lies practically at the surface and is the closest to the training example.
The \gls{DE} deviates further from the \gls{FEM} as its distance to the red dot increases, reaching the maximum of \SI{1.66}{\kelvin \per\centi\meter} at ($\inputCrystalRadius = \SI{76.25}{\milli\meter}$, $\inputReflectorRadius = \SI{39.125}{\milli\meter}$).
On the other hand, \gls{DE} and \gls{FEM} show considerably better agreement across all of $\outputRadialThermalGradient(\mathbf{B})$.
Similarly, $\outputExceedStress(\mathbf{A})$ and $\outputExceedStress(\mathbf{B})$ present excellent agreement between \gls{DE} and \gls{FEM} results.
As observed in \cite{Vieira2025}, $\outputExceedStress$ depends primarily on $\inputCrystalRadius$ and $\inputPullingRate$.

\subsection{Genetic Algorithm Solutions}
\label{sec:GAsolutions}
\begin{figure}
    \centering
    \includegraphics[width=0.85\linewidth]{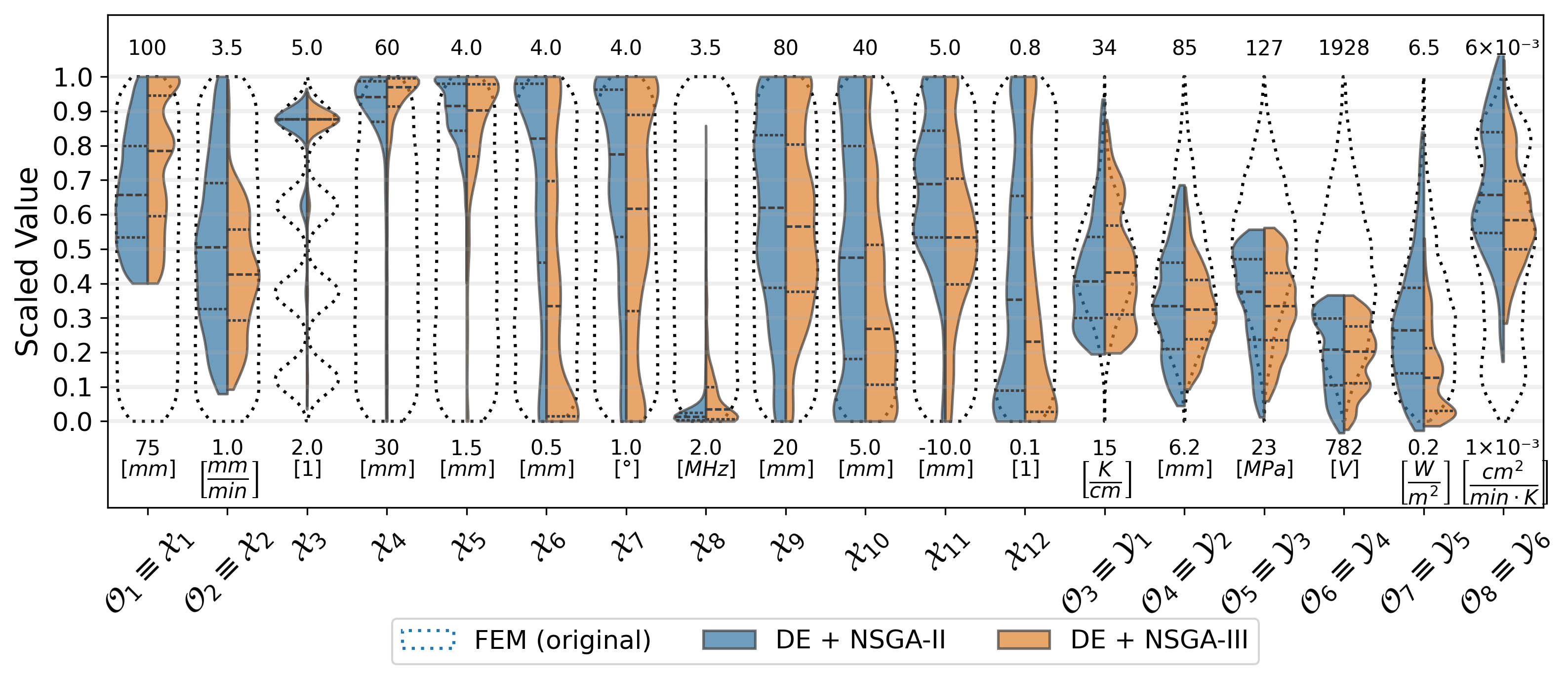}
    \caption{Violin plots with the input and output distributions and their respective ranges.
             Dotted-lines: \gls{FEM} dataset obtained with the \gls{DoE} method.
             Filled regions: optimal solutions obtained with \gls{NSGA}+\gls{DE}.}
    \label{fig:violinplots}
\end{figure}

Figure~\ref{fig:violinplots} overlays the \gls{FEM} results with distributions from the \gls{NSGA}-II and \gls{NSGA}-III optimized solutions (500 solutions each).
All inputs were constrained to the ranges in Table~\ref{tab:modelInputs} and Table~\ref{tab:objectives}.

Clearly, \gls{NSGA}-II and \gls{NSGA}-III agree on $\inputSideSlits \to 5$, $\inputSideSlitLength \to \SI{60}{\milli\meter}$, $\inputSideSlitWidth \to \SI{4}{\milli\meter}$ and $\inputFrequency \to \SI{2.0}{\mega\hertz}$.
The beneficial effect of lowering the working frequency $\inputFrequency$ to reduce the risk of arcing (lower $\outputVoltage$) was investigated in \cite{Menzel2011}, where the authors concluded that this is a viable solution, but further optimization was required.
Increasing the number and length of side slits, $\inputSideSlits$ and $\inputSideSlitLength$, respectively, improves the homogeneity of \gls{EM} field $\outputEMhomogeneity$.
Reducing the width of the main slit $\inputMainSlitWidth$ also improves $\outputEMhomogeneity$, but increases $\outputVoltage$ \cite{Menzel2013}.
\gls{NSGA}-II and \gls{NSGA}-III explored different venues to compensate for this trade-off, as $\inputMainSlitWidth \to \SI{4}{\milli\meter}$ in the former, and $\inputMainSlitWidth \to \SI{0.5}{\milli\meter}$ in the latter.

Objectives $\objective_1 \equiv \inputCrystalRadius$ and $\objective_2 \equiv \inputPullingRate$ reflect process productivity.
$\inputCrystalRadius$ distribution from \gls{NSGA}-III has two well defined modes around 0.8 (\SI{95}{\milli\meter}) and 1.0 (\SI{100}{\milli\meter}), which were traded-off for $\inputPullingRate <$ 0.7 (\SI{2.88}{\milli\meter\per\minute}).
In \gls{NSGA}-II, on the other hand, $\inputCrystalRadius$ and $\inputPullingRate$ distributions are more spreaded out across the whole allowed range.
The elongated tail in $\outputVoronkovCriterion$ from \gls{NSGA}-II, going beyond 1.0 (\SI{5.9e3}{\centi\meter\squared \per\minute\per\kelvin}), stems mainly from the high values of $\inputPullingRate$.

Regarding the outputs, while the constraints in Table~\ref{tab:objectives} help guiding the optimization process and ensure a better balance of objectives, they also prevented the \gls{NSGA} solutions to reach many of the best \gls{FEM} solutions.
In fact, only $\objective_6$ and $\objective_7$ had a slight improvement.
These results indicate a narrow window for process improvement in \gls{FZ}.

The best \gls{NSGA}-III solutions never outperformed their \gls{NSGA}-II counterparts. 
However, \gls{NSGA}-III occasionally found more solutions closer to preset goals (e.g. Q3 and Q1 quartiles in $\objective_1$ and $\objective_7$, respectively).
More reference points in \gls{NSGA}-III should generate ever better solutions.
This is done by increasing the hyperparameter $g$ in Equation~\ref{eq:refPoints}.
Unfortunately, the combinatorial explosion quickly leads to prohibitively long computation times.
We used $p=12$, which generates \num{50388} reference points and corresponds to a 25-fold slowdown compared to \gls{NSGA}-II.

The fact that both \glspl{GA} found solutions already investigated in the past demonstrates the soundness of our approach.
Some nontrivial solutions are discussed and validated in the next Sections.

\subsubsection{Pareto-Optimal Solutions}
\label{sec:optimumSolutions}
\begin{figure}
    \centering
    \includegraphics[width=0.65\linewidth]{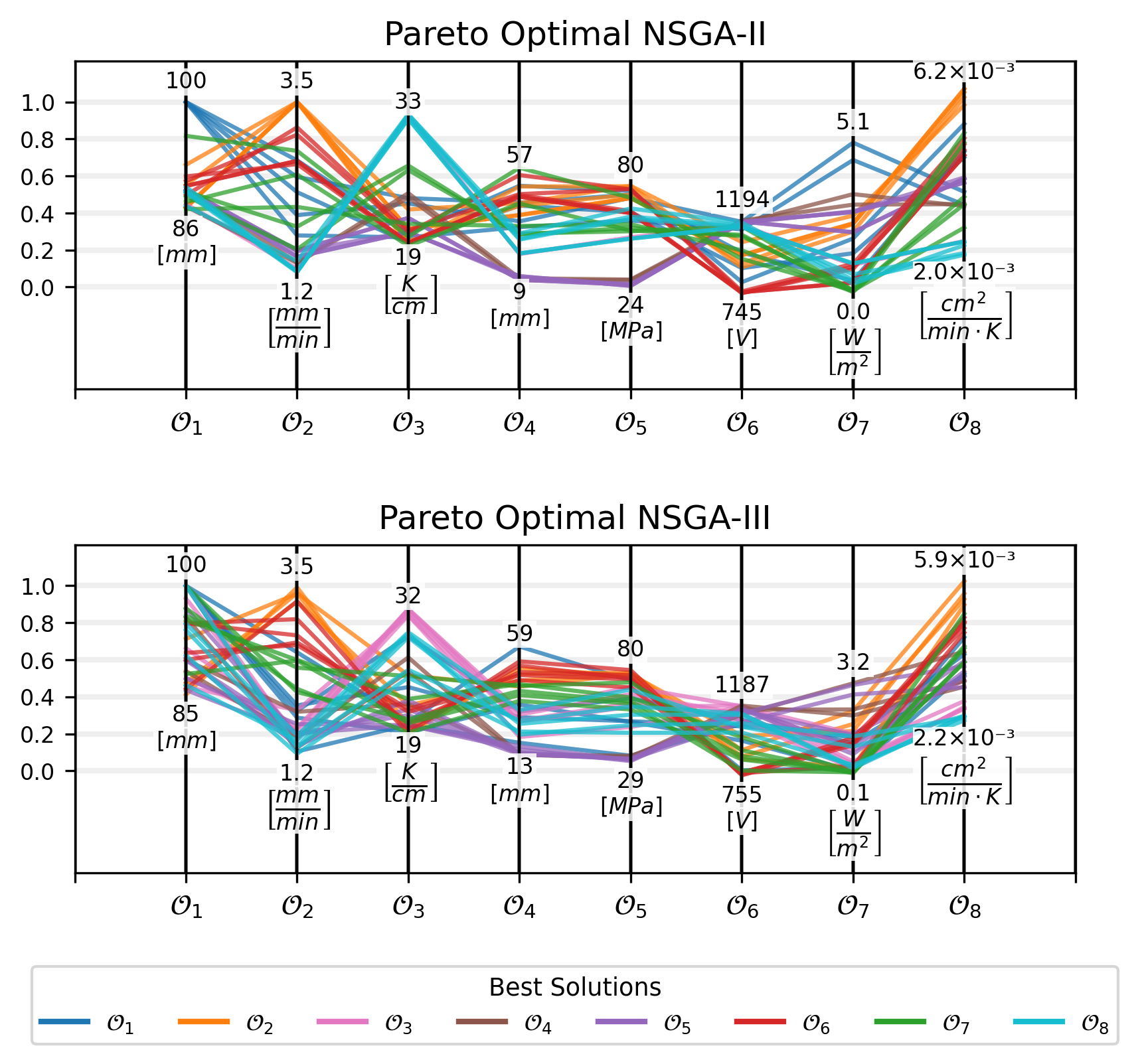}
    \caption{Pareto Optimal solutions obtained with \gls{NSGA}-II and \gls{NSGA}-III. Only five of the best solutions from each objectives are shown.}
    \label{fig:pareto}
\end{figure}

All 500 solutions from \gls{NSGA}-II and \gls{NSGA}-III lie on the first Pareto front, each representing a trade-off among the objectives.
To facilitate visualization, Figure~\ref{fig:pareto} highlights only five of the best solutions for each objective, hiding solutions that exhibit a better balance among all objectives simultaneously.
Some meaningful trade-offs are discussed below.

Tracing the blue lines that emerge from $\inputCrystalRadius \approx \SI{100}{\milli\meter}$ reveals the trade-offs of objective $\objective_1$.
They typically correspond to a pulling rate $\inputPullingRate$ away from the absolute maximum.
Conversely, the highest $\inputPullingRate$ correspond to a low $\inputCrystalRadius$.
The largest Radial Thermal Gradients $\objective_3 \equiv \outputRadialThermalGradient$ in \gls{NSGA}-II have both low $\inputCrystalRadius$ and $\inputPullingRate$, whereas in \gls{NSGA}-III some large values of $\inputCrystalRadius$ are possible.

Exceed stress $\outputExceedStress \equiv \objective_5$ values lower than \SI{30}{\mega\pascal} can only be obtained with small $\inputCrystalRadius$ and $\inputPullingRate$, and they correspond to a small $\outputRadialThermalGradient$ and large voltage $\outputVoltage \equiv \objective_6$.
The deflection $\outputDeflection \equiv \objective_4$ correlates very well to $\outputExceedStress$, as seen by the horizontal lines between $\objective_4$ and $\objective_5$.
Low $\outputVoltage$ are typically associated with low $\outputEMhomogeneity \equiv \objective_7$, which is desirable, but also present low $\outputRadialThermalGradient$ and high $\outputExceedStress$, which are undesirable.
Finally, only \gls{NSGA}-II found solutions with Voronkov's criterion $\outputVoronkovCriterion$ within the critical interval interval \SI{1.3e-3}{} $\le \VoronkovCriterion_\text{crit} \le$ \SI{2.2e-3}{} \SI{}{\square\centi\meter \per\minute \per\kelvin}.
These correspond to very low $\inputPullingRate$, which can induce dislocations \cite{Menzel2013}, and thus should be avoided.
Therefore, it is probably not possible to reliably grow large \gls{FZ} crystals without any kind of point defects.

Figure~\ref{fig:pareto} confirms that \gls{NSGA}-III was able to find slightly more diverse solutions.
However, since \gls{NSGA}-II's complexity scales with the square of population size, comparable results could be obtained by increasing the number of individuals in \gls{NSGA}-II.
Consistent with \cite{Hisao2016}, \gls{NSGA}-III appears unjustifiable for this problem.

\subsubsection{Validation of the Genetic Algorithm Solutions}
\label{sec:validationNSGA2}

\begin{figure}
    \centering
    \includegraphics[width=0.7\linewidth]{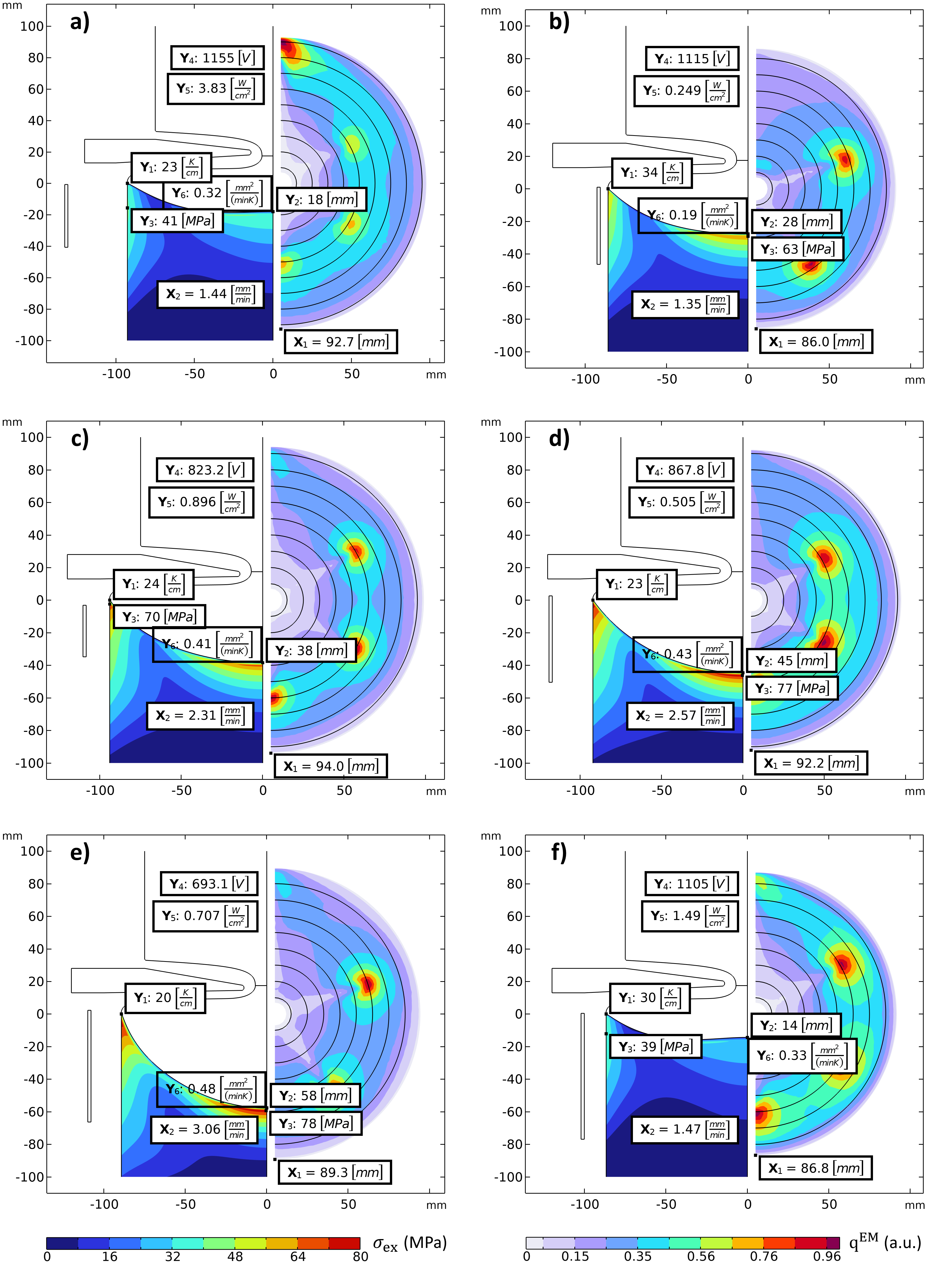}
    \caption{Values of the optimized objectives obtained using the \gls{FEM} for a selection of candidate solutions from \gls{NSGA}-II (subfigures a–f).
    Each subfigure presents simulation results under different parameter settings (e.g., $\inputCrystalRadius$ to $\inputReflectorEmissivity$), with details provided in Table~\ref{tab:validation}.}
    \label{fig:FEM_NSGA_validation}
\end{figure}

To assess the validity of the predictions from the \gls{DE}+\gls{NSGA} framework, some selected cases were recomputed by the \gls{FEM} model.
Figure~\ref{fig:FEM_NSGA_validation} shows the stress field in the crystal and the induced \gls{EM} power on the free surface for six cases.
Due to the symmetries, only half of the fields are shown.
The input values and respective outputs, both computed and predicted, as well as their difference, are detailed in Table~\ref{tab:validation}.

\begin{table}
    \centering
    \begin{tabular}{cccccccc}
     \hline
     & \textbf{\ref{fig:FEM_NSGA_validation}a}) & \textbf{\ref{fig:FEM_NSGA_validation}b}) & \textbf{\ref{fig:FEM_NSGA_validation}c}) & \textbf{\ref{fig:FEM_NSGA_validation}d}) & \textbf{\ref{fig:FEM_NSGA_validation}e}) & \textbf{\ref{fig:FEM_NSGA_validation}f}) & Unit\\
     \hline
    $\inputCrystalRadius$ & 92.7 & 86.0 & 94.0 & 92.2 & 89.3 & 86.8 & \footnotesize\SI{}{\milli\meter} \\
    $\inputPullingRate$ & 1.44 & 1.35 & 2.31 & 2.57 & 3.06 & 1.47 & \footnotesize\SI{}{\milli\meter\per\minute} \\
    $\inputSideSlits$ & 5 & 4 & 5 & 5 & 4 & 5 & \footnotesize\SI{1}{} \\
    $\inputSideSlitLength$ & 49.0 & 55.8 & 58.4 & 50.0 & 58.1 & 59.2 & \footnotesize\SI{}{\milli\meter} \\
    $\inputSideSlitWidth$ & 3.64 & 3.61 & 3.95 & 3.90 & 4.00 & 3.58 & \footnotesize\SI{}{\milli\meter} \\
    $\inputMainSlitWidth$ & 3.44 & 0.55 & 1.89 & 0.81 & 2.04 & 2.13 & \footnotesize\SI{}{\milli\meter} \\
    $\inputBottomAngle$ & 3.40 & 1.59 & 1.00 & 2.31 & 1.00 & 2.79 & \footnotesize\SI{}{\degree}\\
    $\inputFrequency$ & 2.05 & 2.25& 2.08 & 2.01 & 2.01 & 2.00 & \footnotesize\SI{}{\mega\hertz} \\
    $\inputReflectorHeight$ & 39.9 & 47.3 & 31.7 & 53.0 & 68.6 & 77.3 & \footnotesize\SI{}{\milli\meter} \\
    $\inputReflectorRadius$ & 39.9 & 6.78 & 16.4 & 27.0 & 20.7 & 15.4 & \footnotesize\SI{}{\milli\meter} \\
    $\inputReflectorPosition$ & 0.75 & -0.96 & 3.2 & -2.6 & -2.3 & -0.39 & \footnotesize\SI{}{\milli\meter} \\
    $\inputReflectorEmissivity$ & 0.19 & 0.85 & 0.23 & 0.36 & 0.18 & 0.66 & \footnotesize\SI{1}{} \\
    \hline
    $\outputRadialThermalGradient$ & 23.32 & 33.9 & 24.10 & 23.10 & 20.3 & 30.38 & \footnotesize\SI{}{\kelvin\per\centi\meter} \\
    $\hat{\outputRadialThermalGradient}$ & 23.30 & 32.8 & 24.11 & 23.33 & 20.7 & 29.35 & \footnotesize\SI{}{\kelvin\per\centi\meter} \\
     & \footnotesize0.10\% & \footnotesize3.25\% & \footnotesize0.04\% & \footnotesize0.96\% & \footnotesize1.88\% & \footnotesize3.41\% & \\
     \hline
    $\outputDeflection$ & 18.02 & 27.7 & 38.4 & 44.6 & 57.6 & 14.4 & \footnotesize\SI{}{\milli\meter} \\
    $\hat{\outputDeflection}$ & 18.34 & 26.8 & 39.0 & 44.6 & 54.0 & 16.9 & \footnotesize\SI{}{\milli\meter} \\
     & \footnotesize1.78\% & \footnotesize3.30\% & \footnotesize1.35\% & \footnotesize0.11\% & \footnotesize6.28\% & \footnotesize17.2\% & \\
     \hline
    $\outputExceedStress$ & 41.49 & 62.8 & 70.0 & 76.7 & 78.3 & 39.2 & \footnotesize\SI{}{\mega\pascal} \\
    $\hat{\outputExceedStress}$ & 41.98 & 62.6 & 70.1 & 77.5 & 77.9 & 38.9 & \footnotesize\SI{}{\mega\pascal} \\
     & \footnotesize1.08\% & \footnotesize0.43\% & \footnotesize0.08\% & \footnotesize1.04\% & \footnotesize0.51\% & \footnotesize0.79\% & \\
     \hline
    $\outputVoltage$ & 1155 & 1115 & 823.2 & 867.8 & 693.1 & 1105 & \footnotesize\SI{}{\volt} \\
    $\hat{\outputVoltage}$ & 1169 & 1164 & 846.8 & 867.0 & 744.7 & 1147 & \footnotesize\SI{}{\volt} \\
     & \footnotesize1.17\% & \footnotesize4.35\% & \footnotesize2.87\% & \footnotesize0.14\% & \footnotesize7.45\% & \footnotesize3.84\% & \\
     \hline
    $\outputEMhomogeneity$ & 3.83 & 0.249 & 0.896 & 0.505 & 0.707 & 1.49 & \footnotesize\SI{}{\watt\per\centi\meter\squared} \\
    $\hat{\outputEMhomogeneity}$ & 3.80 & 0.280 & 0.934 & 0.484 & 0.787 & 1.51 & \footnotesize\SI{}{\watt\per\centi\meter\squared} \\
     & \footnotesize0.82\% & \footnotesize12.58\% & \footnotesize4.32\% & \footnotesize4.20\% & \footnotesize11.4\% & \footnotesize1.54\% & \\
     \hline
    $\outputVoronkovCriterion$ & 3.177 & 1.865 & 4.075 & 4.326 & 4.833 & 3.326 & \footnotesize \makecell[c]{$\times$1000 \\ \SI{}{\centi\meter\squared \per\minute \per\kelvin}} \\
    $\hat{\outputVoronkovCriterion}$ & 2.963 & 2.192 & 3.966 & 4.386 & 4.823 & 2.948 &  \footnotesize \makecell[c]{$\times$1000 \\ \SI{}{\centi\meter\squared \per\minute \per\kelvin}} \\
     & \footnotesize6.74\% & \footnotesize17.58\% & \footnotesize2.66\% & \footnotesize1.38\% & \footnotesize0.21\% & \footnotesize11.35\% & \\
    \hline
    \end{tabular}
    \caption{Inputs and outputs for each case presented in Figure~\ref{fig:FEM_NSGA_validation}.
    Regarding the outputs, predicted, computed and discrepancy values are provided.}
    \label{tab:validation}
\end{table}

The impact of the main slit is visible at the rim of the free surface in Figure~\ref{fig:FEM_NSGA_validation}a), and measured by $\outputEMhomogeneity = \SI{3.83}{\watt\per\centi\meter\squared}$.
This likely a consequence of the large values for bottom angle and width of the main slit, $\inputBottomAngle$ and $\inputMainSlitWidth$, respectively.
The largest $\inputMainSlitWidth$ was not able to prevent the voltage $\outputVoltage$ from being above $\SI{1100}{\volt}$.
Nevertheless, exceed stress $\outputExceedStress = \SI{41}{\mega\pascal}$ is deemed an acceptable value for a crystal with diameter $2 \times \inputCrystalRadius = \SI{185.5}{\milli\meter}$.

In Figure~\ref{fig:FEM_NSGA_validation}b), while $\outputEMhomogeneity$, the Voronkov's criterion $\outputVoronkovCriterion$ and the Radial Thermal Gradient $\outputRadialThermalGradient$ exhibit the best values for all cases, $\inputCrystalRadius$ and the pulling rate $\inputPullingRate$ present the worst.
Moreover, the reflector is too close to the crystal (offset $\inputReflectorRadius = \SI{6.8}{\milli\meter}$), which poses practical challenges.

Simultaneously reducing $\outputVoltage$, increasing both $\inputCrystalRadius$ and $\inputPullingRate$, and maintaining $\outputEMhomogeneity$ at a relatively low level, was achieved in Figure~\ref{fig:FEM_NSGA_validation}c) and Figure~\ref{fig:FEM_NSGA_validation}d).
However, $\outputExceedStress$, $\outputRadialThermalGradient$ and $\outputVoronkovCriterion$ worsened significantly.
Further increasing $\inputPullingRate$ leads to even worse values for these quantities in Figure~\ref{fig:FEM_NSGA_validation}e), but $\outputVoltage$ exhibits the lowest value.
Finally, Figure~\ref{fig:FEM_NSGA_validation}f) focuses on improving $\outputExceedStress$ and $\outputRadialThermalGradient$, at the expense of inferior $\inputCrystalRadius$, $\inputPullingRate$, $\outputVoltage$ and $\outputEMhomogeneity$. 

The highest absolute discrepancy between predicted and computed outputs, $\hat{\modelOutput}$ and $\modelOutput$, respectively, was \SI{17.58}{\percent} for $\outputVoronkovCriterion$ in case \textbf{\ref{fig:FEM_NSGA_validation}b)}, and the average discrepancy for all cases was \SI{3.84}{\percent}.
However, a significant proportion of the \gls{NSGA} solutions led to nonphysical results when recomputed by the \gls{FEM}.
The reason is twofold: the \gls{DE} was not explicitly trained on the feasibility of the input combinations, and  the \gls{NSGA} was not informed about nonphysical constraints.
Specifically, only \SI{4.4}{\percent} from the initial 2500 sampled combinations yields nonphysical results.
This is not sufficient to effectively train \glspl{NN} on solution feasibility.
Thus, it is imperative that the subset of interest from the \gls{NSGA} solutions be validated against \gls{FEM} simulations before experimental validation.
As an outlook, if a large and diverse set of feasible and nonphysical solutions can be systematically collected during optimization, this information could feature as an additional output in a new \gls{DE}, enabling early prediction of solution validity.

\section{Conclusion}
To support the optimization of Floating Zone (FZ) growth while balancing multiple competing objectives and numerous geometric and growth parameters, we employed a surrogate-based optimization strategy.
To the best of the authors’ knowledge, this is the first systematic investigation of Many-Objective Optimization (MOO) applied to the FZ process.
The surrogate model was built using a Deep Ensemble (DE), a more robust alternative to conventional Neural Networks, trained on numerical data from a Finite Element Model (FEM) of the FZ process.
MOO was carried out using two variants of the Non-dominated Sorting Genetic Algorithm (NSGA-II and NSGA-III), generating a set of candidate solutions, each optimal for one or more objectives.

Between the two optimization strategies, NSGA-II outperformed NSGA-III, despite not being specifically designed for high-dimensional problems.
The candidate solutions aligned with findings from previous studies, such as the reduction of working frequency, supporting the validity of our approach.
The Pareto frontier revealed both expected trade-offs (e.g. crystal size versus thermal stress), as well as more intricate ones.
A subset of candidate solutions was validated by recomputing the outputs using FEM simulations.
Overall, the agreement between the DE and FEM was excellent, with an average discrepancy of less than \SI{4}{\percent}.
Nevertheless, because solution feasibility was not explicitly included as a constraint, selected candidate solutions must be re-evaluated with the FEM before moving to experimental validation.

This framework provides a robust and efficient tool for identifying promising parameter combinations in crystal growth optimization for FZ, and can be extended to other growth methods.
However, our results indicate that the process improvement window for FZ is narrow.
Further improvements will require incorporating additional information into the surrogate model.


\medskip
\textbf{Data Availability} \par Data will be made available on request.

\medskip
\textbf{Conflict of Interest} \par The authors declare that no conflicts of interest exist.

\bibliographystyle{MSP}
\bibliography{bibliography}

\medskip
\textbf{Appendix}

\par \textbf{A.}

\begin{table}[H]
    \centering
    \sisetup{
          scientific-notation = true,  
          round-mode = figures,       
          round-precision = 3         
    }
    \begin{tabular}{cccccccccccc}
        \hline
         \multirow{2}{*}{\textbf{MSE}} & \multirow{2}{*}{\textbf{Hidden}} & \multicolumn{10}{c}{\textbf{Neurons per Hidden Layer}} \\
         \cmidrule(ll){3-12}
         ($\times$ 10\textsuperscript{-4}) & \textbf{Layers} & \textbf{1} & \textbf{2} & \textbf{3} & \textbf{4} & \textbf{5} & \textbf{6} & \textbf{7} & \textbf{8} & \textbf{9} & \textbf{10}\\
         \hline
         \num{3.4432} & 6 & 62 & 63 & 60 & 58 & 61 & 54 &  &  &  &  \\
         \num{3.4499} & 5 & 56 & 53 & 64 & 54 & 63 &  &  &  &  & \\
         \num{3.7879} & 4 & 41 & 56 & 60 & 64 &  &  &  &  &  & \\
         \num{4.0297} & 7 & 47 & 61 & 38 & 59 & 64 & 43 & 48 &  &  & \\
         \num{4.1025} & 3 & 37 & 62 & 63 &  &  &  &  &  &  & \\
         \num{4.5247} & 10 & 37 & 62 & 63 & 59 & 58 & 62 & 47 & 53 & 46 & 37\\
         \num{4.8282} & 8 & 40 & 53 & 58 & 56 & 59 & 58 & 57 & 51 &  & \\
         \num{5.0644} & 2 & 50 & 61 &  &  &  &  &  &  &  & \\
         \num{5.9809} & 9 & 54 & 14 & 60 & 54 & 58 & 39 & 14 & 27 & 47 & \\
         \num{7.3511} & 1 & 37 &  &  &  &  &  &  &  &  & \\
         \hline
    \end{tabular}
    \caption{Best-performing NN architectures (selected from 1000 candidates) for \glspl{NN} with 1 to 10 hidden layers.
Noticeably, architectures containing latent representations (e.g. few neurons in a given layer) is unfavorable in our setting.}
    \label{tab:bestNNs}
\end{table}

\par \textbf{B.}
\begin{table}[H]
    \centering
    \begin{tabular}{cccc}
        \hline
         & \textbf{A} & \textbf{B} & \textbf{Unit} \\
        \hline
        $\inputCrystalRadius$ & 95.3* & 77.2* & \SI{}{\milli\meter}\\
        $\inputPullingRate$ & 1.1* & 2.3* & \SI{}{\milli\meter \per\minute}\\
        $\inputSideSlits$ & 5 & 2  & \SI{}{1}\\
        $\inputSideSlitLength$ & 45.4 & 43.8 & \SI{}{\milli\meter}\\
        $\inputSideSlitWidth$ & 2.9 & 2.2 & \SI{}{\milli\meter}\\
        $\inputMainSlitWidth$ & 3.8* & 1.7* & \SI{}{\milli\meter}\\
        $\inputBottomAngle$ & 3.9 & 3.2 &  \SI{}{\degree}\\
        $\inputFrequency$ & 2.7 & 2.5 & \SI{}{\mega\hertz}\\
        $\inputReflectorHeight$ & 57.5 & 26.8 & \SI{}{\milli\meter}\\
        $\inputReflectorRadius$ & 25.2* & 38.4* & \SI{}{\milli\meter}\\
        $\inputReflectorPosition$ & 8.9 & 5.8 & \SI{}{\milli\meter}\\
        $\inputReflectorEmissivity$ & 0.4 & 0.8 & \SI{}{1}\\
        \hline
    \end{tabular}
    \caption{Input combinations (\textbf{A}, \textbf{B} $\in \trainingData$) used to compute the \gls{DE} as a function of two inputs. Inputs with * have their value overwritten by the respective range in Figure~\ref{fig:NN_surfacePlots}.}
    \label{tab:combinationsSurfacePlot}
\end{table}

\begin{figure}
\textbf{Table of Contents}\\
\centering
\medskip
  \includegraphics[width=0.85\linewidth]{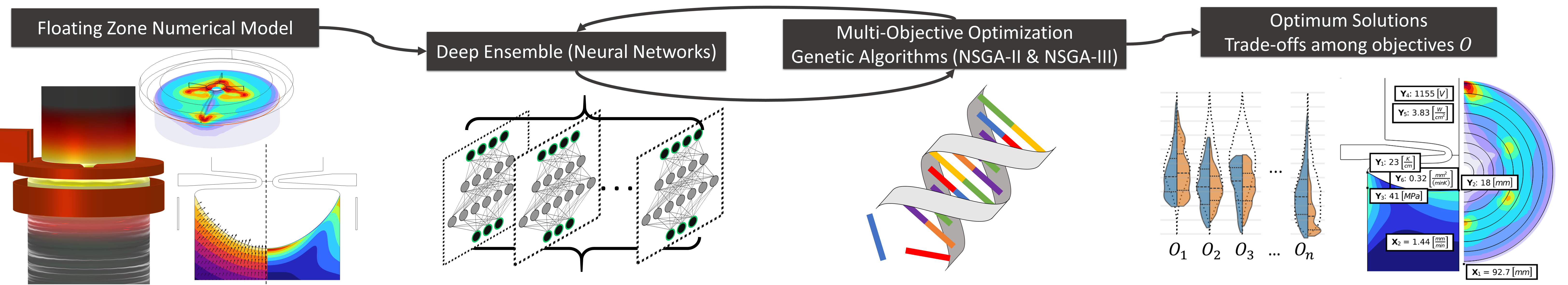}
  \medskip
  \caption*{This study presents a surrogate-based Multi-Objective Optimization framework for Floating Zone silicon crystal growth. An ensemble of Neural Networks is trained on simulation data and combined with Genetic Algorithms to explore trade-offs in process parameters. Results highlight the predictive power of the model and its potential for guiding crystal growth process improvements.}
\end{figure}

\end{document}